\documentclass[aps, showpacs, showkeys, onecolumn, amsmath, amssymb, superscriptaddress, preprintnumbers, floatfix, preprint]{revtex4}

\usepackage{bm,amsfonts, mathtools}
\usepackage{graphicx,color, wasysym}
\usepackage{textcomp}
\usepackage{hyperref} 
\usepackage{msc}

\def\de{\mathrm{d}}
\def\e{\mathrm{e}}
\def\nn{\nonumber}

\numberwithin{equation}{section}

\begin{document}

\title{Repeated derivatives of composite functions and generalizations of the Leibniz rule}

\author{\large{D. Babusci}}
\email{danilo.babusci@lnf.infn.it}
\affiliation{INFN - Laboratori Nazionali di Frascati, via E. Fermi, 40, I 00044 Frascati (Roma), Italy}

\author{\large{G. Dattoli}}
\email{dattoli@frascati.enea.it}
\affiliation{ENEA - Centro Ricerche Frascati, via E. Fermi, 45, I 00044 Frascati (Roma), Italy}

\author{\large{K. G\'{o}rska}}
\email{kasia\_gorska@o2.pl}
\affiliation{Instituto de F\'{\i}sica, Universidade de S\~{a}o Paulo, P.O.Box 66318, BR 05315-970 S\~{a}o Paulo, SP, Brasil}
\affiliation{H. Niewodnicza\'{n}ski Institute of Nuclear Physics, Polish Academy of Sciences, ul.Eljasza-Radzikowskiego 152, 
PL 31342 Krak\'{o}w, Poland}
\affiliation{Laboratoire de Physique Th\'{e}orique de la Mati\`{e}re Condens\'{e}e (LPTMC), Universit\'{e} Pierre et Marie Curie, 
CNRS UMR 7600, Tour 13 - 5i\`{e}me \'{e}t., Bo\^{i}te Courrier 121, 4 place Jussieu, F 75252 Paris Cedex 05, France}

\author{\large{K. A. Penson}}
\email{penson@lptl.jussieu.fr}
\affiliation{Laboratoire de Physique Th\'{e}orique de la Mati\`{e}re Condens\'{e}e (LPTMC), Universit\'{e} Pierre et Marie Curie, 
CNRS UMR 7600, Tour 13 - 5i\`{e}me \'{e}t., Bo\^{i}te Courrier 121, 4 place Jussieu, F 75252 Paris Cedex 05, France}

\begin{abstract}
We use the properties of Hermite and Kamp\'e de F\'eriet polynomials to get closed forms for the repeated derivatives of 
functions whose argument is a quadratic or higher-order polynomial. The results we obtain are extended to product of functions 
of the above argument, thus giving rise to expressions which can formally be interpreted as generalizations of the familiar Leibniz 
rule. Finally, examples of practical interest are discussed.
\end{abstract}
\maketitle

\section{Introduction}
Formula (1.1.1.1) of Ref. \cite{YuABrychkov08} refers to repeated derivatives of functions whose argument is a quadratic polynomial, 
i.e. ($\hat{D}_\xi = \de/\de \xi$)
\begin{equation}\label{eqDer}
\hat{D}_x^n\,[f(x^{2})] = n!\,\sum_{k = 0}^{[n/2]} \frac{(2 x)^{n - 2 k}}{k!\, (n - 2 k)!}\, \hat{D}_{x^2}^{n-k}\,f (x^2).
\end{equation}
Albeit elementary, and a particular case of the Fa\`{a} di Bruno formula \cite{FdiBruno}, some aspects of Eq. (\ref{eqDer}) 
deserve to be studied as they may lead to novel results.

Eq. (\ref{eqDer}) implicitly assumes that $f(x)$ is at least a $n$-times differentiable function. Moreover, we make the 
hypothesis that it can be written as
\begin{equation}\label{eqf}
f(x) = \sum_{n = 0}^{\infty} \frac{x^n}{n!}\,f_n = \e^{x\,\hat{\phi}}\,f_0, \qquad\qquad \hat{\phi}^{n} f_{0} = f_{n},
\end{equation}
where $f_{0}$ denotes a kind of ``vacuum" state on which the repeated action of the umbral operator $\hat{\phi}$ generates 
a discrete sequence of functions. In what follows we will consider only functions that admit the expansion (\ref{eqf}), i.e., that are 
analytical over the entire complex plane.

As an example, we note that the Tricomi function of order zero \cite{Kasia} can be written as
\begin{equation}
C_{0} (x) = \sum_{k = 0}^{\infty} \frac{(-x)^k}{(k!)^2} = \e^{- x\,\hat{\phi}}\,f_0, 
\end{equation}
with
\begin{equation}
\hat{\phi}^k\,f_{0} = \frac1{\Gamma (k + 1)}, \qquad\qquad f_{0} = \frac1{\Gamma(1)} = 1.
\end{equation}
The operator $\hat{\phi}$ may also be raised to any real (not necessarily integer) power, so that the Tricomi function 
of order $\alpha$ can be written as
\begin{equation}
C_{\alpha} (x) = \sum_{k = 0}^{\infty} \frac{(-x)^k}{k!\,\Gamma(k + \alpha + 1)} = \hat{\phi}^\alpha \,\e^{- x\,\hat{\phi}}\,f_0\,,
\end{equation}
The obvious advantage of the previous umbral re-shaping of the function $f (x)$ is the possibility of exploiting the wealth of 
the properties of the exponential function, which will be assumed to be still valid, since we will treat the operator $\hat{\phi}$ 
as an ordinary constant \cite{Kasia}.

As a consequence of Eq. (\ref{eqf}), we can write (we omit the vacuum $f_0$ for brevity)
\begin{equation}\label{eqDgx}
\hat{D}_{g (x)}^n\,f [g (x)] = \hat{\phi}^{n}\,\e^{\,g(x)\,\hat{\phi}}.
\end{equation}
By taking into account the identity \cite{Appell}
\begin{equation}\label{eqHer}
\hat{D}_x^{n}\,\e^{\,a\,x^2} =  H_{n}^{(2)} (2\,a\,x, a)\,\e^{\,a\,x^2}, 
\end{equation}
where
\begin{equation}
H_n^{(2)} (x, y) = n! \sum_{k = 0}^{[n/2]} \frac{x^{n - 2 k}\, y^{k}}{(n - 2 k)!\, k!}
\end{equation}
is the two-variable Hermite-Kamp\'e de F\'eriet polynomials \cite{Appell, BDDbook},  one obtains  
\begin{eqnarray}\label{eqDxn}
\hat{D}_x^n\,f(x^2) = \hat{D}_x^n\,\e^{\,x^2\,\hat{\phi}} &=& H_n^{(2)} (2\,x\,\hat{\phi}, \hat{\phi})\,
\e^{x^2\,\hat{\phi}} \nn \\ 
&=& \left\{n!\,\sum_{k = 0}^{[n/2]} \frac{(2 x)^{n - 2 k}}{(n - 2 k)!\, k!}\, \hat{\phi}^{n - k}\right\}\,\e^{\,x^2\,\hat{\phi}}\,,
\end{eqnarray}
from which, by using Eq. (\ref{eqDgx}), we recover Eq. (\ref{eqDer}) (see also Appendix A).

\section{An extension of Leibniz formula}
In this section we will draw further consequences from the previous formalism. In particular, we will consider the case in which 
the function $f (x)$ is the product of two functions, obtaining a closed formula which will be recognized as a generalized version 
of the Leibniz rule.

We start by considering the following function
\begin{equation}
f (x^2)= g (x^2)\,h (x^2) = \left(\e^{\,x^2\,\hat{\gamma}}\, g_0\right)\, \left(\e^{\,x^2\,\hat{\eta}}\,h_0\right),
\end{equation}
i.e., since the operators $\hat{\gamma}$ and $\hat{\eta}$ commute (we omit the vacuum terms)
\begin{equation}
f (x^2) = \e^{\,x^2\,(\hat{\gamma} + \hat{\eta})}\,,
\end{equation}
and, according to the Eq. (\ref{eqDxn}), we find
\begin{eqnarray}\label{eqDngh}
\hat{D}_x^n\, [f (x^2)] &=& \hat{D}_x^{n}\,\e^{\,x^2\,(\hat{\gamma} + \hat{\eta})} = H_{n}^{(2)}\left(2(\hat{\gamma} + \hat{\eta})\,x, 
\hat{\gamma} + \hat{\eta}\right)\,\e^{\,x^2\,(\hat{\gamma} + \hat{\eta})} \nn \\
&=& \left\{n!\, \sum_{k = 0}^{[n/2]}\,\frac{(2 x)^{n - 2 k}}{(n - 2 k)!\, k!}\, (\hat{\gamma} + \hat{\eta})^{n - k}\,\right\}\,
\e^{\,x^2\,(\hat{\gamma} + \hat{\eta})} \\
&=& \left\{n! \sum_{k = 0}^{[n/2]} \frac{(2 x)^{n - 2 k}}{(n - 2 k)!\, k!}\, \sum_{j = 0}^{n - k} \binom{n - k}{j}\, \hat{\gamma}^{n - k- j}\, 
\hat{\eta}^j\right\}\,\e^{\,x^2\,(\hat{\gamma} + \hat{\eta})}\,. \nn
\end{eqnarray}
By noting that
\begin{equation}
\hat{\gamma}^{\,m}\,\hat{\eta}^{\,p}\,\e^{\,x^2\,(\hat{\gamma} + \hat{\eta})} = (\hat{\gamma}^{\,m}\,\e^{\,x^2\,\hat{\gamma}}\, g_0)\,
(\hat{\eta}^{\,p}\,\e^{\,x^2\,\hat{\eta}}\,h_0) = \left[\hat{D}_{x^2}^m\,g (x^2)\right]\,\left[\hat{D}_{x^2}^p\,h (x^2)\right],
\end{equation}
we obtain
\begin{equation}
\hat{D}_x^{n}\,[g (x^2)\,h (x^2)] = n! \sum_{k = 0}^{[n/2]} \frac{(2 x)^{n - 2 k}}{(n - 2 k)!\,k!}\,\sum_{j = 0}^{n - k} \binom{n - k}{j}\, 
\left[\hat{D}_{x^2}^{n - k - j}\,g (x^2)\right]\,\left[\hat{D}_{x^2}^j\,h (x^2)\right],
\end{equation}
which is manifestly a generalization of the ordinary Leibniz rule for the $n$th  derivative of the product of two functions. The 
use of the addition formula \cite{BDDbook}
\begin{equation}
H_{n}^{(2)}(x_1 + x_2, y_1 + y_2) = \sum_{k = 0}^{n} \binom{n}{k}\, H_{n - k}^{(2)}(x_1, y_1)\,H_{k}^{(2)}(x_2, y_2)
\end{equation}
in Eq. (\ref{eqDngh}) yields
\begin{equation}
\hat{D}_x^{n}\,[g (x^2)\,h (x^2)] = \sum_{k = 0}^{n} \binom{n}{k}\,\left[H_{n - k}^{(2)}(2\,x\,\hat{\gamma}, \hat{\gamma})\,
\e^{\,x^2\,\hat{\gamma}}\right]\,\left[H_k^{(2)}(2\,x\,\hat{\eta}, \hat{\eta})\,\e^{\,x^2\,\hat{\eta}}\right],
\end{equation}
which makes the analogy more transparent. From this example, we conclude that, for the type of function we are dealing with, 
the operator 
\begin{equation}\label{eq2Gam}
_{(2)}\hat{\Gamma}^{(n)} = H_n^{(2)} (2\,\hat{\gamma}\,x, \hat{\gamma})
\end{equation}
is a kind of multiple derivative operator. For more detail on this point see Appendix B. 

A fairly interesting application of Eq. (\ref{eqDngh}) is in the calculation of $n$th derivative of the product of cylindrical 
Bessel functions, that, in our present formalism, writes \cite{Kasia}
\begin{equation}\label{eqBes}
J_{n} (x) = \left(\hat{\varphi}\,\frac{x}2\right)^n\,\e^{- (x/2)^2\,\hat{\varphi}}\,j_{0}, \qquad\qquad 
\hat{\varphi}^{n}\, j_0 = j_{n} = \frac1{\Gamma(n + 1)}.
\end{equation}
One obtains
\begin{equation}
\hat{D}_x^{n}\,[J^2_0 (x)] = (- 1)^n\,n!\, \sum_{k = 0}^{[n/2]} \frac{(- 2 x)^{- k}}{(n - 2 k)!\, k!}\, 
\sum_{m = 0}^{n - k} \binom{n - k}{m}\,J_{n - k - m} (x)\,J_{m} (x).
\end{equation}
Obviously, the formalism can easily applied also to other products of functions.

We close this section, by considering $n$-th derivative of the function $f(a x^2 + b x)$. From Eq. (\ref{eqf}), one has
\begin{equation}
\hat{D}_x^n\,[f (a x^2 +  b x)] = \hat{D}_x^n\,\e^{\,(a x^2 + b x)\,\hat{\chi}}\,f_0
\end{equation}
and the use of the identity
\begin{equation}
\hat{D}_x^{n}\,\e^{a x^2 + b x} = H_n^{(2)}(2\,a\,x + b, a)\,\e^{a x^2 + b x}
\end{equation}
leads to
\begin{equation}
\hat{D}_x^{n} [f(a x^2 + b x)] = n! \sum_{k = 0}^{[n/2]} \frac{(2 a x + b)^{n - 2 k}\, a^{k}}{(n - 2 k)!\,k!}\,
\hat{D}_{a x^2 + b x}^{n - k}(a x^2 + b x),
\end{equation}

\section{Higher-order Hermite polynomials and repeated derivatives}
We have already remarked that the examples we are considering here are the particular cases of the Fa\`{a} di Bruno formula, which 
had different formulations in the course of the last two centuries \cite{FdiBruno}. Among the various possibilities there is that of 
expressing the $n$th derivative of a composite functions in terms of the Bell polynomials \cite{Bell}, which are a generalization of 
the Hermite-Kamp\'e de F\'eriet polynomials. A family of polynomials intermediate between that exploited so far and 
the Bell one is represented by the three variable Hermite-Kamp\'e de F\'eriet polynomials, defined as 
\cite{Pino99}
\begin{equation}
H_{n}^{(3)} (x_1, x_2, x_3) = n!\,\sum_{k = 0}^{[n/3]} \frac{x_3^k\,H_{n - 3 k}^{(2)}(x_1, x_2)}{(n - 3 k)!\,k!}
\end{equation}
The following identity generalizes Eq. (\ref{eqHer}) (see Refs. \cite{BDDbook} and \cite{Pino99})
\begin{equation}
\hat{D}_x^n\,\e^{a x^3} = H_n^{(3)} (3\,a\,x^2, 3\,a\,x, a)\,\e^{a x^3}.
\end{equation}
and, accordingly, it is evident that
\begin{equation}\label{eqDfx3}
\hat{D}_x^n f (x^3) = H_n^{(3)} (3\,x^2\,\hat{\phi}, 3\,x\,\hat{\phi}, \hat{\phi})\,\e^{\,x^3\,\hat{\phi}}, 
\end{equation}
i.e.,
\begin{equation}
\hat{D}_x^{n} f (x^{3}) = \left\{n!\,\sum_{k = 0}^{[n/3]}\,\sum_{m = 0}^{[(n - 3 k)/2]} 
\frac{(3 x^2)^{n - 3 k - m}\,x^{- m}}{(n - 3 k - 2 m)!\,k!\,m!}\,\hat{D}_{x^3}^{n - 2 k - m}\right\}\,f (x^3).
\end{equation}
It is evident that the method we are developing is a kind of modular scheme, which can be easily generalized and automatized. 
An example is provided by the following Leibniz rule
\begin{equation}
\hat{D}_x^{n} [g (x^3)\,h (x^3)] = \sum_{k = 0}^{n} \binom{n}{k}\,H_{n - k}^{(3)} (3\,x^2\,\hat{\gamma}, 3\,x\,\hat{\gamma}, \hat{\gamma})\,
H_k^{(3)} (3\,x^2\,\hat{\eta}, 3\,x\,\hat{\eta}, \hat{\eta}), 
\end{equation}
where the addition formula
\begin{equation}
H_n^{(3)}(x_1 + x_2, y_1 + y_2, z_1 + z_2) = \sum_{k = 0}^{n} \binom{n}{k}\,H_{n - k}^{(3)} (x_1, y_1, z_1)\, 
H_k^{(3)} (x_2, y_2, z_2)
\end{equation}
has been used.

\section{Concluding Remarks}
In this section we complete the previous discussion by using the properties of the higher-order Hermite polynomials. 
In addition, we will comment on other formulae reported in Ref. \cite{YuABrychkov08} about successive derivatives of functions.

By using the identity \cite{BDDbook,Pino99,Pino03}
\begin{equation}
\hat{D}_x^{n}\,\e^{\,a\,x^m} = H_n^{(m)} (a\,_m z_1, \ldots, a\,_m z_m)\,\e^{\,a\,x^m} \qquad\qquad _m z_k = \binom{m}{k}\,x^{m - k},
\end{equation}
with
\begin{equation}
H_n^{(m)} (\xi_1, \ldots, \xi_m) = n!\, \sum_{k = 0}^{[n/m]} \frac{\xi_m^k}{(n - m k)!\,k!}\,H_{n - m k}^{(m - 1)}(\xi_1, \ldots, \xi_{m - 1}), 
\end{equation}
we obtain the following generalization of Eq. (\ref{eqDfx3}): 
\begin{eqnarray}
\hat{D}_x^n\ f (x^m) &=& n!\,\sum_{k_1 = 0}^{[n/m]} \frac{_m z_m^{\,k_1}}{k_1!}\,\sum_{k_2 = 0}^{[\ell_{m - 1}]} 
\frac{_m z_{m - 1}^{\,k_2}}{k_2!}\,\ldots\,\sum_{k_{m - 2} = 0}^{[\ell_3]}\frac{_m z_3^{\,k_{m - 2}}}{k_{m - 2}!}\,
\sum_{k_{m - 1} = 0}^{[\ell_2]} \frac{_m z_2^{\,k_{m - 1}}}{k_{m - 1}!}\,\frac{_m z_1^{\,\ell_1}}{\ell_1!}\,
\hat{D}_{x^m}^{\,\ell_1 + p_1}\,f (x^m) \nn \\
&& \qquad \qquad \left(\ell_j = \frac1{j}\,[n - \sum_{p = j}^{m - 1} (p + 1)\,k_{m - p}], \quad 
p_1 = \sum_{q = 0}^{m - 1} k_{m - q}\right).
\end{eqnarray}

As stated in the incipit, this paper has been inspired by Eq. (1.1.1.1) of Ref. \cite{YuABrychkov08}. We pass now to consider  
Eq. (1.1.1.2), which states that, for $n \geq 1$,
\begin{equation}\label{eqrx}          
\hat{D}_x^n\,f(\sqrt{x}) = \sum_{k = 0}^{n-1}\,(-1)^k\,\frac{(n + k - 1)!}{k!\,(n - k - 1)!}\,\frac1{(2\,\sqrt{x})^{n + k}}\,
\hat{D}_{\sqrt{x}}^{n - k} f(\sqrt{x}).
\end{equation}
We can write
\begin{eqnarray}
\sum_{n = 0}^{\infty} \frac{t^n}{n!}\,\hat{D}_x^{n}\,f(\sqrt{x}) &=& \e^{\,t\,\hat{D}_x}\,f (\sqrt{x}) = f(\sqrt{x + t}) = 
\e^{\,\sqrt{x + t}\,\hat{\phi}}\,f_0 \nn \\
&=& \e^{\,\sqrt{x}\,\hat{\phi}}\,\exp\left\{- \sqrt{x}\,\left(1 - \sqrt{1 +  \frac{t}{x}}\right)\,\hat{\phi}\right\}\,f_0, 
\end{eqnarray}
and, taking into account the expression of the generating function for the Bessel polynomials (see formula (25) of Ref. \cite{HLKrall45})
\begin{equation}                                                                     
\sum_{n = 0}^{\infty}\,\frac{t^n}{n!}\,y_{n - 1} (x) = \exp\left\{\frac1{x}\,(1 - \sqrt{1 - 2 x t})\right\}, \qquad 
y_n (x) = \sum_{k = 0}^{n}\,\frac{(n + k)!}{(n - k)!\,k!}\,\left(\frac{x}2\right)^k,
\end{equation}
we obtain 
\begin{equation}
\sum_{n = 0}^{\infty}\,\frac{t^n}{n!}\,\hat{D}_x\,f (\sqrt{x}) = \sum_{n = 0}^{\infty}\,\frac1{n!}\,
\left(\frac{t\,\hat{\phi}}{2\,\sqrt{x}}\right)^n\,y_{n - 1}\left(\frac{-1}{\sqrt{x}\,\hat{\phi}}\right)\,\e^{\,\sqrt{x}\,\hat{\phi}}\,f_0.
\end{equation}
from which Eq. (\ref{eqrx}) is easily obtained equating the coefficients of the same powers in $t$.

As a final example of application of the method, we prove formula (1.1.1.3) of Ref. \cite{YuABrychkov08},
\begin{equation}\label{eqf1x}
\hat{D}_x^n\,f (1/x) = (- 1)^n\,(n - 1)!\,\sum_{k = 0}^n \binom{n}{k}\,\frac{x^{k - 2 n}}{(n - k -1)!}\,\hat{D}_{1/x}^{n - k}\,f (1/x)\,, 
\qquad (n \leq 1)\,.
\end{equation}
Indeed, after setting
\begin{equation}
\xi = \frac1{x}, \qquad\qquad f (\xi) = \e^{\,\hat{c}\,\xi}, \nn
\end{equation} 
we can write
\begin{eqnarray}\label{eqgen}
\sum_{n = 0}^\infty \frac{t^{\,n}}{n!}\,\hat{D}_x^n\,f (1/x) &=&  \sum_{n = 0}^\infty (- 1)^n\,\frac{t^{\,n}}{n!}\,
(\xi^2\,\hat{D}_\xi)^n\, f (\xi) = \e^{- t\,\xi^2\,\hat{D}_\xi}\,\e^{\,\hat{c}\,\xi} \nn \\
&=& \exp\left\{\frac{\xi\,\hat{c}}{1 + \xi\,t}\right\} = \exp\left\{- \frac{\xi^2\,t\,\hat{c}}{1 + \xi\,t}\right\}\,\e^{\,\xi\,\hat{c}}, 
\end{eqnarray}
where in the second line we have used the identity \cite{Comtet}
\begin{equation}
\e^{\,\lambda\,x^2\,\hat{D}_x}\,f (x) = f \left(\frac{x}{1 - \lambda\,x}\right)\,, \qquad\qquad \left(|\lambda\,x| < 1\right).
\end{equation}
By remembering that \cite{Pino99}
\begin{equation}
\sum_{n = 0}^\infty t^{\,n}\,L_n (x,y) = \frac1{1 - y\,t}\,\exp\left\{- \frac{x\,t}{1 - y\,t}\right\},
\end{equation}
where
\begin{equation}
L_n (x, y) = n!\,\sum_{k = 0}^n (- 1)^k\,\frac{x^{\,k}\,y^{\,n - k}}{(n - k)!\,(k!)^2} 
\end{equation}
are the two-variable Laguerre polynomials, it's easy to show that (with $L_{-1} (x, y) = 0$)
\begin{eqnarray}
\exp\left\{- \frac{\xi^2\,t\,\hat{c}}{1 + \xi\,t}\right\} &=& \sum_{n = 0}^\infty t^n\,\left[L_n (\xi^{\,2}\,\hat{c}, - \xi) + 
\xi\,L_{n - 1} (\xi^{\,2}\,\hat{c}\,, - \xi)\right] \nn \\
&=& 1 + \sum_{n = 1}^\infty t^n\,\frac{(- 1)^n}n\,\sum_{k = 0}^n \binom{n}{k}\,\frac{\xi^{\,2 n - k}}{(n - k - 1)!}\,\hat{c}^{\,n - k} 
\end{eqnarray}
from which, going back to the variable $x$, taking into account Eq. (\ref{eqDgx}), and comparing powers of $t$ of the same order in (\ref{eqgen}), 
Eq. (\ref{eqf1x}) is obtained.


\appendix

\section{}
By setting $x^2 = y$ in Eq. (\ref{eqDer}) we get the following operatorial identity
\begin{equation}
\left(\sqrt{y}\,\hat{D}_y\right)^n = \sum_{k = 0}^{[n/2]} S_2^{(1/2)} (n, k)\,(\sqrt{y})^{n - 2 k}\,\hat{D}_y^{n - k}
\end{equation}
where
\begin{equation}
S_2^{(1/2)} (n, k) = \frac{n!}{4^k\,k!\,(n - 2 k)!}.
\end{equation}

Analogously, by setting $\sqrt{x} = z$ in Eq. (\ref{eqrx}) we get ($n \geq 1$)
\begin{equation}
\left(\frac1{z}\,\hat{D}_z\right)^n = \sum_{k = 0}^{n - 1} S_2^{(- 1)} (n, k)\,z^{- n - k}\,\hat{D}_z^{n - k}
\end{equation}
where
\begin{equation}
S_2^{(- 1)} (n, k) = (- 1)^k\,\frac{(n + k - 1)!}{2^k\,k!\,(n - k - 1)!}.
\end{equation}

The numbers $S_2^{(\nu)} (n, k)$ are a generalization of the Stirling numbers of second kind $S_2 ^{(1)} (n, k)$ \cite{Comtet} involved 
in the expansion
\begin{equation}
\left(x\,\hat{D}_x\right)^n = \sum_{k = 0}^{n} S_2^{(1)} (n, k)\,x^k\,\hat{D}_x^k.
\end{equation}

\section{}
It is well known that symbols in mathematics have their own life and meaning. We can therefore generalize the operator 
defined in Eq. (\ref{eq2Gam}) as follows
\begin{equation}
_{(2)}\hat{\Gamma}^{(n)} (a, b) = H_n^{(2)} (2\,a\,x\,\hat{\gamma}, b\,\hat{\gamma})
\end{equation}
and note that
\begin{equation}
_{(2)}\hat{\Gamma}^{(n)} (a/2, b)\,f (x^2) = \left\{n!\,\sum_{k = 0}^{[n/2]} \frac{a^{n - 2 k}\,b^k}{(n - 2 k)!\,k!}\, \hat{D}_{x^2}^{n - k}\right\}\,f (x^2).
\end{equation}
If the function $f$ can be expanded as in Eq. \eqref{eqf}, we can consider the problem of evaluating the following integral
\begin{eqnarray}
\mathcal{I}_n (a, b, c) &=& \int_{-\infty}^{\infty}\de x\,_{(2)}\hat{\Gamma}^{(n)} (a/2, b)\,f(- c\,x^2)  \nn \\
&=& \left(\int_{-\infty}^{\infty}\de x\,H_n^{(2)} (a\,x\,\hat{\gamma}, b\,\hat{\gamma})\,\e^{- c\,x^2\,\hat{\gamma}}\right)\,f_0.
\end{eqnarray}
By remembering that the generating function of the two-variable Hermite-Kamp\'e de F\'eriet polynomials is given by
\begin{equation}
\sum_{n = 0}^{\infty} \frac{t^n}{n!}\,H_n^{(2)} (y, z) = \e^{\,y\,t + z\,t^2},
\end{equation}
we get
\begin{equation}
\sum_{n = 0}^{\infty} \frac{t^n}{n!}\,\mathcal{I}_n (a, b, c) = \left(\int_{-\infty}^{\infty} \de x\,\e^{\,a\,x\,t\,\hat{\gamma} + b\,t^2\,\hat{\gamma}}\, 
\e^{- c\,x^2\,\hat{\gamma}}\right)\,f_0
\end{equation}
i.e., treating the integral on the r. h. s. as an ordinary Gaussian integral
\begin{equation}
\sum_{n = 0}^{\infty} \frac{t^{n}}{n!}\,\mathcal{I}_{n}(a, b, c) = \sqrt{\frac{\pi}{c\,\hat{\gamma}}}\,
\exp\left[\left(\frac{a^2}{4\,c} + b\right)\,t^2\,\hat{\gamma}\right]\,f_0.
\end{equation}
In this equation, by expanding the exponential in series and equating the $t$-like power terms, we find ($n \leq 1$)
\begin{equation}
\mathcal{I}_{2 n} (a, b, c) = 2^n\,(2 n - 1)!!\,\sqrt{\frac{\pi}{c}}\,\left(\frac{a^2}{4\,c} + b\right)^n\,f_{n - 1/2}, 
\qquad\qquad \mathcal{I}_{2 n + 1} (a, b, c) = 0.
\end{equation}
As an application of this result, from Eq. (\ref{eqBes}) one has
\begin{equation}
\int_{-\infty}^{\infty}\de x\,_{(2)}\hat{\Gamma}^{(2 n)} (a/2, b)\,J_0 (x) = 2^{\,n + 1}\,(2 n - 1)!!\,\sqrt{\pi}\,
\frac{(a^{2} + b)^n}{\Gamma (n + 1/2)}. 
\end{equation}

\vspace{0.5cm}
\begin{center}
\textbf{Acknowledgements}\\
\end{center}
\noindent
G. D. dedicates this paper to his friend Paolo Iudici, whose longstanding friendship was one of the firm points in his life. 
K. G. thanks Funda\c{c}\~{a}o de Amparo \'{a} Pesquisa do Estado de S\~{a}o Paulo (FAPESP, Brazil) under Program 
No. 2010/15698-5. K. G. and K. A. P. acknowledge support from Agence Nationale de la Recherche (Paris, France) under 
Program PHYSCOMB No. ANR-08-BLAN-243-2.

\end{document}